# On computing differential transform of nonlinear non-autonomous functions and its applications

Essam. R. El-Zahar[1,2] and Abdelhalim Ebaid[3]

[1]Department of Mathematics, Faculty of Sciences and Humanities, Prince Sattam Bin Abdulaziz University, Alkharj, 11942, KSA.

[2]Department of Basic Engineering Science, Faculty of Engineering, Shebin El-Kom, 32511, Menofia University, Egypt.

[3]Department of Mathematics, Faculty of Science, Tabuk University, P. O. Box 741, Tabuk 71491, KSA

**Abstract**: Although being powerful, the differential transform method yet suffers from a drawback which is how to compute the differential transform of nonlinear non-autonomous functions that can limit its applicability. In order to overcome this defect, we introduce in this paper a new general formula and its related recurrence relations for computing the differential transform of any analytic nonlinear non-autonomous function with one or multi-variable. Regarding, the formula in the literature was found not applicable to deal with the present non-autonomous functions. Accordingly, a generalization is presented in this paper which reduces to the corresponding formula in the literature as a special case. Several test examples for different types of nonlinear differential and integro-differential equations are solved to demonstrate the validity and applicability of the present method. The obtained results declare that the suggested approach not only effective but also a straight forward even in solving differential and integro-differential equations with complex nonlinearities.

## 1. Introduction

The Differential Transform Method (DTM) which is based on Taylor series expansion was first introduced by Zhou [1] and has been successfully applied to a wide class of nonlinear problems arising in mathematical sciences and engineering. The main advantage of the DTM is that it can be applied directly to nonlinear differential equations with no need for linearization, discretization, or perturbation. Additionally, the DTM does not generate secular terms (noise terms) and does not need to analytical integrations as other semi-analytical numerical methods such as HPM, HAM, ADM or VIM and so the DTM is an attractive tool for solving differential equations. Although this method has been proved to be an efficient tool for handling nonlinear differential equations, the nonlinear functions used in these studies are restricted to certain kinds of nonlinearities, e.g., polynomials and products with derivatives. For other types of nonlinear functions, Chang and Chang [2] construct a new algorithm based on obtaining a differential equation satisfied by this nonlinear function and then applying the DTM to this obtained differential equation. Although their treatment was found effective for some forms of nonlinearity [3,4] it significantly increases the computational budget, especially if there are two or more nonlinear functions involved in the differential equation being investigated, as demonstrated by Ebaid [5]. Moreover, in the case of complex nonlinearities, it may be quite difficult to obtain the differential equations satisfied by these nonlinear functions. To overcome this difficulty, a new formula has been derived [5-8] to calculate the differential transform of any nonlinear autonomous one variable function $f(y)$. This formula has the same mathematical structure as the Adomian polynomials but with constants instead of variable components. We mean by the nonlinear non-autonomous functions that those nonlinear functions which contain the independent variable and also the dependent variable and/or its derivatives. Unfortunately, for such types of complicated nonlinear non-autonomous functions with one or multi-variable $f(t, y_j(t))$, $j = 1, 2, ..., m$, no related formula has been given to calculate their transform functions. This provides the motivation for the present work.

To overcome the drawback of the DTM that is how to compute the differential transform of high complex nonlinear non-autonomous functions, a generalized formula is deduced in this paper for computing the differential transform of any analytic nonlinear non-autonomous function with one or multi-variable. The proposed method deals directly with the nonlinear non-autonomous function in its form without any special kinds of transformations or algebraic manipulations. Also, there is no need to compute the differential transform of other functions to obtain the required one. For autonomous function, as a special case of the current study, these formulas and recurrence relations have the same mathematical structure as the Adomian polynomials but with constants instead of variable components. The efficiency of the proposed method is discussed through several test examples including nonlinear differential and integro-differential equations of different types.

## 2. Differential Transform Method

The basic definitions and fundamental theorems of the one-dimensional DTM and its applicability for various kinds of differential and integro-differential equations are given in [1, 9-15]. For convenience of the reader, we present a brief review of the DTM in this section.

The differential transform of a given analytic function $y(t)$ in a domain $D \subset \mathbb{R}$ is defined as

$$Y(i) = \left[\frac{1}{i!}\frac{d^i y(t)}{dt^i}\right]_{t=t_0}, \quad t_0 \in D, \quad (1)$$

where $y(t)$ is the original function and $Y(i)$ is the transformed function. The inverse differential transform of $Y(i)$ is defined as

$$y(t) = \sum_{i=0}^{\infty} Y(i)(t-t_0)^i. \quad (2)$$

From Eqs. (1) and (2), we get

$$y(t) = \sum_{i=0}^{\infty} \left[\frac{1}{i!}\frac{d^i y(t)}{dt^i}\right]_{t=t_0} (t-t_0)^i. \quad (3)$$

From the above proposition, it can be found that the concept of differential transform is derived from Taylor series expansion. In actual applications, the function $y(t)$ is expressed by a truncated series and Eq. (2) can be written as

$$\tilde{y}(t) = \sum_{i=0}^{N} Y(i)(t-t_0)^i, \quad (4)$$

where $N$ is the approximation order of the solution. Some of the fundamental mathematical operations performed by the one dimensional differential transform are listed in Table 1.

**Table 1.** Some fundamental operations of the one dimensional differential transform.

| Original function $y(t)$ | Transformed function $Y(i)$ |
|---|---|
| $\beta(v(t) \pm w(t))$ | $\beta V(k) \pm \beta W(k)$ |
| $v(t)w(t)$ | $\sum_{\ell=0}^{i} V(\ell)W(i-\ell)$ |
| $\dfrac{d^m v(t)}{dt^m}$ | $\dfrac{(i+m)!}{i!}V(i+m)$ |
| $(\beta+t)^m$ | $H[m,i]\dfrac{m!}{i!(m-i)!}(\beta+t_0)^{m-i}$, where $H[m,i]=\begin{cases} 1, & \text{if } m \geq i \\ 0, & \text{if } m < i \end{cases}$ |
| $e^{\lambda t}$ | $\dfrac{\lambda^k}{i!}e^{\lambda t_0}$ |
| $\int_{t_0}^{t} v(\tau)d\tau$ | $\dfrac{V(i-1)}{i}, i \geq 1$ |
| $u(t)\int_{t_0}^{t} v(\tau)d\tau$ | $\sum_{\ell=0}^{i-1} G(\ell)\dfrac{V(i-\ell-1)}{i-\ell}, i \geq 1$ |
| $\dfrac{d^m v(qt)}{dt^m}$ | $\dfrac{(i+m)!}{i!}q^{i+m}V(i+m)$ |
| $\sin(\omega t + \beta)$ | $\dfrac{\omega^i}{i!}\sin(\omega t_0 + \beta + \dfrac{i\pi}{2})$ |
| $\cos(\omega t + \beta)$ | $\dfrac{\omega^i}{i!}\cos(\omega t_0 + \beta + \dfrac{i\pi}{2})$ |

## 3. Differential Transform formulas

If a differential equation contains an analytic nonlinear non-autonomous function $f(t, y(t))$ then the differential transform $F(n)$ of the function $f(t, y(t))$ is required and can be computed from the following theorems, where we assume that $f(t, y(t)) \neq f_0(t)f_1(y(t))$.

**Theorem 1**. *The differential transform $F(n)$ of any analytic nonlinear non-autonomous function $f(t, y(t))$ at a point $t_0$ can be computed from the formula*

$$F(n) = \frac{1}{n!}\frac{d^n}{d\lambda^n}\left[ f\left(t_0 + \lambda, \sum_{i=0}^{\infty} Y(i)\lambda^i\right)\right]_{\lambda=0} \quad (4)$$

where $Y(i)$ is the differential transform of $y(t)$.

**Proof**:

The differential transform $F(n)$ of $f(t, y(t))$ is defined as

$$F(n) = \frac{1}{n!}\left[\frac{d^n}{dt^n}f(t, y(t))\right]_{t=t_0}. \quad (5)$$

Substituting Eq.2 in Eq.5 results in

$$F(n) = \frac{1}{n!}\left[\frac{d^n}{dt^n}f\left(t,\sum_{i=0}^{\infty}Y(i)(t-t_0)^i\right)\right]_{t=t_0}. \tag{6}$$

Now, let $t - t_0 = \lambda$, then Eq.6 becomes

$$F(n) = \frac{1}{n!}\left[\frac{d^n}{d\lambda^n}f\left(t_0+\lambda,\sum_{i=0}^{\infty}Y(i)\lambda^i\right)\right]_{\lambda=0}.$$

By this way, the proof of Theorem 1 is completed.

**Theorem 2.** *The differential transform $F(n)$ of any analytic nonlinear non-autonomous function $f(t,y(t))$ at a point $t_0$, satisfies the recurrence relation*

$$F(n) = \frac{1}{n}\left(\frac{\partial}{\partial t_0}F(n-1) + \sum_{i=0}^{n-1}(i+1)Y(i)\frac{\partial}{\partial Y(i)}F(n-1)\right), n=1,2,..., \tag{7}$$

where $F(0) = f(t_0, Y(0))$.

**Proof:** we have
$$f^{(n)}(t,y(t)) = \frac{\partial}{\partial t}f^{(n-1)}(t,y(t)) + \frac{\partial}{\partial y(t)}f^{(n-1)}(t,y(t))\frac{dy(t)}{dt}, \text{ then}$$

$$F(n) = \frac{1}{n!}\left[\frac{\partial}{\partial t}f^{(n-1)}\left(t,\sum_{i=0}^{\infty}Y(i)(t-t_0)^i\right) + \sum_{i=0}^{\infty}\frac{\partial}{\partial Y(i)}f^{(n-1)}\left(t,\sum_{i=0}^{\infty}Y(i)(t-t_0)^i\right)(i+1)Y(i+1)\right]_{t=t_0}$$

$$= \frac{1}{n!}\left[(n-1)!\frac{\partial}{\partial t_0}F(n-1) + (n-1)!\sum_{i=0}^{\infty}(i+1)Y(i)\frac{\partial}{\partial Y(i)}F(n-1)\right]$$

and since $F(n-1)$ is a function of $t_0$ and $\{Y(i)\}_{i=0}^{n-1}$, then

$$F(n) = \frac{1}{n}\left(\frac{\partial}{\partial t_0}F(n-1) + \sum_{i=0}^{n-1}(i+1)Y(i)\frac{\partial}{\partial Y(i)}F(n-1)\right)$$

By this way, the proof of Theorem 2 is completed.

Thus by Theorems 1 and 2 we have implemented a new algorithm for computing the one-dimensional differential transform of any analytic nonlinear function non-autonomous $f(t,y(t))$.

We observe that for autonomous function, the definition (4) and recurrence relation (7) are reduced to, respectively

$$F(n) = \frac{1}{n!}\frac{d^n}{d\lambda^n}\left[f\left(\sum_{i=0}^{\infty}Y(i)\lambda^i\right)\right]_{\lambda=0}, \tag{8}$$

$$F(n) = \frac{1}{n}\left(\sum_{k=0}^{n-1}(k+1)Y(k)\frac{\partial}{\partial Y(k)}F(n-1)\right), n=1,2,... \tag{9}$$

One can observe that the result in (8) is the present formula in [5-8]

In fact if a system of one-dimensional differential equations contains a coupled analytic nonlinearity $f(t, y_j(t))$, $j = 1,2,...,m$, where we assume that $f(t, y_j(t)) \neq f_0(t) f_1(y_1(t)) f_2(y_2(t)) ... f_m(y_m(t))$, then Theorems 1 and 2 can be easily extended using the above present procedure to multi-variable function and satisfy the following recurrence algorithms

**Corollary 1.** *The differential transform $F(n)$ of any analytic non-autonomous function $f(t, y_j(t))$, $j = 1,2,...,m$ at a point $t_0$ can be defined by*

$$F(n) = \frac{1}{n!} \frac{d^n}{d\lambda^n} \left[ f\left(t_0 + \lambda, \sum_{i=0}^{\infty} Y_j(i) \lambda^i \right) \right]_{\lambda=0}, \quad j = 1,2,...,m \tag{10}$$

where $Y_j(i)$ is the differential transform of $y_j(t)$.

**Corollary 2.** *The differential transform $F(n)$ of any analytic non-autonomous function $f(t, y_j(t))$, $j = 1,2,...,m$ at a point $t_0$, satisfies the recurrence relation*

$$F(n) = \frac{1}{n} \left( \frac{\partial}{\partial t_0} F(n-1) + \sum_{j=1}^{m} \sum_{k=0}^{n-1} (k+1) Y_j(k) \frac{\partial}{\partial Y_j(k)} F(n-1) \right), \quad n = 1,2,..., \tag{11}$$

where $F(0) = f(t_0, Y_j(0))$.

Moreover for autonomous function, the definition (10) and recurrence relation (11) are reduced, respectively, to

$$F(n) = \frac{1}{n!} \frac{d^n}{d\lambda^n} \left[ f\left( \sum_{i=0}^{\infty} Y_j(i) \lambda^i \right) \right]_{\lambda=0}, \quad j = 1,2,...,m, \tag{12}$$

$$F(n) = \frac{1}{n} \left( \sum_{j=1}^{m} \sum_{k=0}^{n-1} (k+1) Y_j(k) \frac{\partial}{\partial Y_j(k)} F(n-1) \right), \quad j = 1,2,...,m, \ n = 1,2,..., \tag{13}$$

which have the same mathematical structure as the Adomian polynomials [16] but with constants instead of variable components.

## 4. Results

In this section, we have solved different types of differential and integro-differential problems with different forms of nonlinear non-autonomous functions to demonstrate the validity and applicability of the present method.

**Example 1.** Consider the nonlinear initial-value problem

$$y'(t) - y(t) = \ln(t + y(t)), \quad y(1) = 0. \tag{14}$$

Using the basic properties of DTM and taking the transform of equations in (14) result in

$$(k+1)Y(k+1) - Y(k) = F(k), \quad Y(0) = 0, \quad k = 0,1,2,\ldots, \qquad (15)$$

where $F(k)$ is the differential transform of the nonlinear term $\ln(t + y(t))$. $F(k)$ is computed using the present method and given by

$$\left.\begin{aligned}
F(0) &= \ln(1+Y(0)), \quad F(1) = \frac{1+Y(1)}{1+Y(0)}, \quad F(2) = \frac{Y(2)}{1+Y(0)} - \frac{1}{2}\frac{(1+Y(1))^2}{(1+Y(0))^2}, \\
F(3) &= \frac{Y(3)}{1+Y(0)} - \frac{Y(2)(1+Y(1))}{(1+Y(0))^2} + \frac{(1+Y(1))^3}{3(1+Y(0))^3}, \\
F(4) &= \frac{Y(4)}{1+Y(0)} - \frac{Y(3)(1+Y(1))}{(1+Y(0))^2} + \frac{Y(2)(1+Y(1))^2}{(1+Y(0))^3} - \frac{Y(2)^2}{2(1+Y(0))^2} - \frac{(1+Y(1))^4}{4(1+Y(0))^4} \\
&\vdots
\end{aligned}\right\} \qquad (16)$$

Therefore, a combination of (15) and (16) results in the series solution

$$\tilde{y}(t) = \frac{1}{2}(t-1)^2 + \frac{1}{6}(t-1)^3 + \frac{1}{24}(t-1)^4 + \frac{1}{120}(t-1)^5 + \ldots.$$

For sufficiently large number of terms, the closed form of the solution is $\tilde{y}(t) = e^{t-1} - t$, which is the exact solution. Table 2 shows the absolute error obtained for three various numbers of terms and at some test points.

**Table 2:** Numerical results for Example 1.

| | $|y(t_i) - \tilde{y}(t_i)|$ | | |
|---|---|---|---|
| $t_i$ | $N = 5$ | $N = 10$ | $N = 15$ |
| 1.0 | 0.00000000 | 0.00000000 | 0.00000000 |
| 1.2 | 9.1494e-008 | 6.0021e-016 | 7.9797e-017 |
| 1.4 | 6.0310e-006 | 1.0869e-012 | 4.1633e-017 |
| 1.6 | 7.0800e-005 | 9.5652e-011 | 5.5511e-017 |
| 1.8 | 4.1026e-004 | 2.3048e-009 | 1.4988e-015 |
| 2.0 | 1.6152e-003 | 2.7313e-008 | 5.1181e-014 |

**Example 2**. Consider the nonlinear initial-value problem

$$y'(t) + \varepsilon y(t)^2 = \varepsilon \sin(ty(t)), \qquad y(0) = 1. \qquad (17)$$

Taking the differential transform of equations in (17) results in

$$(k+1)Y(k+1) + \varepsilon \sum_{\ell=0}^{k} Y(\ell)Y(k-\ell)Y(k) = \varepsilon F(k), \quad Y(0) = 1, \quad k = 0,1,2,\ldots, \qquad (18)$$

where $F(k)$ is the differential transform of the nonlinear term $\sin(ty(t))$. $F(k)$ is computed using the present method and given by

$$F(0) = 0, \quad F(1) = Y(0), \quad F(2) = Y(1), \quad F(3) = \frac{-Y(0)^3}{6} + Y(2),$$
$$F(4) = y(3) - \frac{1}{2}Y(0)^2 y(1), \quad F(5) = \frac{1}{120}Y(0)^5 - \frac{1}{2}Y(0)Y(1)^2 - \frac{1}{2}Y(0)^2 Y(2) + Y(4) \right\} \quad . \tag{19}$$
$$\vdots$$

Using (18) and (19), the series solution is obtained and given at $\varepsilon = 0.1$ by

$$\tilde{y}(t) = 1 - \frac{1}{10}t + \frac{3}{50}t^2 - \frac{23}{3000}t^3 - \frac{119}{60000}t^4 + \frac{247}{300000}t^5 - \frac{2233}{4500000}t^6 + \ldots \quad .$$

The presented results are compared with those obtained using MATLAB built-in solver ode45 in Table 3. The ode45 solver integrates ODEs using explicit 4th & 5th Runge-Kutta (4, 5) formula [17]. In order to guarantee a good numerical reference, ode45 is configured using an absolute error of $10^{-12}$ and relative error of $10^{-8}$.

**Table 3:** Numerical results for Example 2.

| | $\|y(t_i) - \tilde{y}(t_i)\|$ | | |
|---|---|---|---|
| $t_i$ | $N = 5$ | $N = 10$ | $N = 15$ |
| 0.0 | 0.00000000 | 0.00000000 | 0.00000000 |
| 0.2 | 7.3689e-006 | 8.2100e-010 | 9.6316e-011 |
| 0.4 | 4.5510e-004 | 1.3793e-006 | 2.3291e-008 |
| 0.6 | 4.8955e-003 | 1.1150e-004 | 3.0230e-006 |
| 0.8 | 2.5675e-002 | 2.4110e-003 | 2.6524e-004 |
| 1.0 | 9.0957e-002 | 2.5516e-002 | 8.3823e-003 |

**Example 3.** Consider the following nonlinear problem with multiple solutions

$$(t+1)y'(t) = \sqrt{t + y(t)^2}, \quad y(1) = 0. \tag{20}$$

The differential transform of equations in (20) are

$$\sum_{\ell=0}^{k}(\ell+1)Y(\ell+1)H(1,k-\ell)\frac{(2)^{1-k+\ell}}{(k-\ell)!|(1-k+\ell)|!} = F(k) \text{ and } Y(0) = 0,$$

where $F(k)$ is the differential transform of the nonlinear term $\sqrt{t + y(t)^2}$. Applying the present method to the nonlinear function $f(t, y(t)) = \sqrt{t + y(t)^2}$ at $t_0 = 1$ and $Y(0) = 0$, results in

$$F(0) = \pm 1, \; F(1) = \pm\frac{1}{2}, \; F(2) = \pm\left(-\frac{1}{8} + \frac{Y(1)^2}{2}\right), \; F(3) = \pm\left(\frac{1}{16} - \frac{1}{4}Y(1)^2 + Y(1)Y(2)\right)$$
$$F(4) = \pm\left(-\frac{1}{8}Y(1)^4 + \frac{3}{16}Y(1)^2 + \frac{1}{128}(128Y(3) - 64Y(2))Y(1) - \frac{5}{128} + \frac{1}{2}Y(2)^2\right),$$
$$F(5) = \pm\left(\frac{3Y(1)^4}{16} - \frac{Y(1)^3 Y(2)}{2} - \frac{5Y(1)^2}{32} + \left(Y(4) - \frac{Y(3)}{2} + \frac{3Y(2)}{8}\right)Y(1) + \frac{7}{256} - \frac{Y(2)^2}{4} + Y(2)Y(3)\right)$$
$$\vdots \right\} \quad . \tag{21}$$

where the brackets with positive sign are due to using the principal square roots while the brackets with negative sign are due to using the negative roots. Therefore, a combination of (20) and (21) results in two series solution given by

$$\tilde{y}(t) = \pm\left(\frac{1}{2}t - \frac{1}{2}\right),$$

which are the exact solutions of (20).

**Example 4.** Consider the nonlinear first order Volterra integro-differential equation

$$y'(t) = \cos t - \frac{t^2}{2} + 1 + \int_0^t \sin^{-1}(1-\tau+y(\tau))d\tau, \qquad y(0)=-1, \ y'(0)=2. \tag{22}$$

The differential transform of equations in (22) are

$$(k+1)Y(k+1) = \frac{1}{k!}\cos\left(\frac{k\pi}{2}\right) - \frac{\delta(k-2)}{2} + \delta(k) + \frac{F(k-1)}{k}, \ Y(0)=-1, Y(1)=2, k\geq 1,$$

where $F(k)$ is the differential transform of the nonlinear term $\sin^{-1}(1-\tau+y(\tau))$. By applying the present method to the nonlinear function $f(\tau, y(\tau)) = \sin^{-1}(1-\tau+y(\tau))$ at $\tau_0 = 0$, $Y(0)=-1$ and $Y(1)=2$, and by using the principal value of the square root and inverse trigonometric functions, we get

$$\left.\begin{array}{l} F(0)=0, F(1)=1, \quad F(2)=Y(2), \ F(3)=Y(3)+\frac{1}{6}, F(4)=Y(4)+\frac{1}{2}Y(2), \\ F(5)=Y(5)+\frac{3}{40}+\frac{1}{2}Y(3)+\frac{1}{2}Y(2)^2, \ F(6)=\frac{3}{8}Y(2)+\frac{1}{2}Y(4)+Y(2)Y(3)+Y(6)+\frac{1}{6}Y(2)^3 \\ F(7)=\frac{1}{112}(56Y(3)+84)Y(2)^2+Y(2)Y(4)+\frac{3}{8}Y(3)+\frac{1}{2}Y(5)+Y(7)+\frac{1}{2}Y(3)^2+\frac{5}{112} \\ \vdots \end{array}\right\}. \tag{23}$$

Hence, the series solution is obtained and given by

$$\tilde{y}(t) = -1 + 2t - \frac{1}{6}t^3 + \frac{1}{120}t^5 - \frac{1}{5040}t^7 + \frac{1}{362880}t^9 + \ldots$$

For sufficiently large number of terms, the closed form of the solution is $\tilde{y}(t) = \sin(t) + t - 1$, which is the exact solution. Table 4 shows the absolute error obtained for three various numbers of terms and at some test points.

**Table 4:** Numerical results for Example 4.

| | $|y(t_i) - \tilde{y}(t_i)|$ | | |
|---|---|---|---|
| $t_i$ | $N=5$ | $N=10$ | $N=15$ |
| 0.0 | 0.00000000 | 0.00000000 | 0.00000000 |
| 0.2 | 2.5383e-009 | 5.5511e-016 | 1.0011e-017 |
| 0.4 | 3.2436e-007 | 1.0497e-012 | 5.5511e-017 |
| 0.6 | 5.5266e-006 | 9.0679e-011 | 1.0002e-017 |
| 0.8 | 4.1242e-005 | 2.1432e-009 | 1.1102e-016 |
| 1.0 | 1.9568e-004 | 2.4892e-008 | 2.7756e-015 |

**Example 5.** Consider the nonlinear second order Volterra integro-differential equation

$$y''(t) - 2y(t)y'(t) = -t + \int_0^t \frac{\sec^2 \tau}{1 + y(\tau)^2} d\tau, \qquad y(0) = 0, \ y'(0) = 1. \qquad (24)$$

The differential transform of equations in (24) are

$$(k+2)(k+1)Y(k+2) = \left(2\sum_{\ell=0}^{k}(\ell+1)Y(\ell+1)Y(k-\ell) - \delta(k-1) + \frac{F(k-1)}{k}\right), \ Y(0)=0, Y(1)=1, Y(2)=0, k \geq 1$$

, where $F(k)$ is the differential transform of the nonlinear term $\frac{\sec^2 \tau}{1 + y(\tau)^2}$ obtained using the present method at $\tau_0 = 0$ and given by

$$\begin{aligned}
&F(0)=1, F(1)=0, \ F(2)=0, \ F(3)=0, \ F(4)=\frac{2}{3}-2Y(3), Y(5)=-2Y(4), \\
&F(6)=\frac{-13}{45}-2Y(5)+2Y(3)-Y(3)^2, F(6)=(2-2Y(3))Y(4)-2Y(6), \\
&F(7)=\frac{17}{35}+5Y(3)^2-\frac{10}{3}Y(3)-2Y(3)Y(5)+2Y(5)-2Y(7)-Y(4)^2 \\
&\vdots
\end{aligned} \qquad (25)$$

Hence, the approximate series solution is obtained and given by

$$\tilde{y}(t) = t + \frac{1}{3}t^3 + \frac{2}{15}t^5 + \frac{17}{315}t^7 + \frac{62}{2835}t^9 + \ldots$$

For sufficiently large number of terms, the closed form of the solution is $\tilde{y}(t) = \tan t$, which is the exact solution. Table 5 shows the absolute error obtained for three various numbers of terms and at some test points.

**Table 5:** Numerical results for Example 5.

| | $\lvert y(t_i) - \tilde{y}(t_i) \rvert$ | | |
|---|---|---|---|
| $t_i$ | $N = 5$ | $N = 10$ | $N = 15$ |
| 0.0 | 0.00000000 | 0.00000000 | 0.00000000 |
| 0.2 | 7.0218e-007 | 1.8451e-010 | 8.0491e-016 |
| 0.4 | 9.4552e-005 | 3.9753e-007 | 1.0839e-010 |
| 0.6 | 1.7688e-003 | 3.7649e-005 | 1.1693e-007 |
| 0.8 | 1.5281e-002 | 1.0280e-003 | 1.7939e-005 |
| 1.0 | 9.0741e-002 | 1.4903e-002 | 9.9212e-004 |

**Example 6.** Consider the nonlinear Volterra integro-differential equation with proportional delay

$$y'\left(\frac{t}{2}\right) = \frac{1}{2} - t \sin t + \int_0^t \frac{y(3\tau)^2 \sin t}{(3\tau+1)^2} d\tau, \qquad y(0)=1, y'(0)=1. \qquad (26)$$

The differential transform of equations in (26) are

$$(k+1)Y(k+1)\left(\frac{1}{2}\right)^{k+1} = \frac{1}{2}\delta(k) - \sum_{\ell=0}^{k}\frac{1}{\ell!}\sin\left(\frac{\ell\pi}{2}\right)\delta(k-\ell-1) + \sum_{\ell=0}^{k-1}\frac{1}{\ell!}\sin\left(\frac{\ell\pi}{2}\right)\frac{F(k-\ell-1)}{k-\ell}, \ Y(0)=1, Y(1)=1, k \geq 1,$$

where $F(k)$ is the differential transform of the nonlinear term $\dfrac{w(\tau)^2}{(3\tau+1)^2}$, $w(\tau) = \sum_{k=0}^{\infty} 3^k Y(k)\tau^k$ obtained using the present method and given by

$$\begin{aligned}
&F(0)=1,\ F(1)=0,\ F(2)=18Y(2),\ F(3)=-54Y(2)+54Y(3),\\
&F(4)=81Y(2)^2+162Y(2)-162Y(3)+162Y(4),\\
&F(5)=-486Y(2)^2+(486Y(3)-486)Y(2)-486Y(4)+486Y(5)+486Y(3)\\
&\vdots
\end{aligned} \qquad (27)$$

Hence, the series solution is obtained and given by

$$\tilde{y}(t) = t+1,$$

which is the exact solution.

**Example 7.** Consider the following nonlinear non-autonomous initial-value ODE system

$$\left.\begin{aligned}
y_1'(t) &= -y_1(t) + t + \ln\!\left(y_1(t) - \dfrac{1}{t+y_2(t)}\right)\\
y_2'(t) &= -y_2(t) - 1 + \dfrac{4}{y_1(t)} - \ln(t+y_2(t))\\
y_1(0) &= 2,\quad y_2(0) = 1
\end{aligned}\right\} . \qquad (28)$$

Applying the differential transform to (28), results in

$$\left.\begin{aligned}
(k+1)Y_1(k+1) &= \delta(K-1) - Y_1(k) + F_1(k)\\
(k+1)Y_2(k+1) &= -Y_2(k) - \delta(K) + F_2(k)\\
Y_1(0) &= 2,\quad Y_2(0) = 1
\end{aligned}\right\} . \qquad (29)$$

where $F_1(k)$ and $F_2(k)$ are the differential transform of the nonlinear functions $f_1 = \ln\!\left(y_1(t) - \dfrac{1}{t+y_2(t)}\right)$ and $f_2 = \dfrac{4}{y_1(t)} - \ln(t+y_2(t))$, respectively. $F_1(k)$ and $F_2(k)$ are computed using formula(10) and given by

$$\left.\begin{aligned}
&F_1(0) = 0,\ F_1(1) = Y_1(1) + 1 + Y_2(1),\\
&F_1(2) = Y_1(2) - (1+Y_2(1))^2 + Y_2(2) - \dfrac{1}{2}(Y_1(1)+1+Y_2(1))^2\\
&F_1(3) = Y_1(3) + (1+Y_2(1))^3 - 2(1+Y_2(1))Y_2(2) + Y_2(3) -\\
&\qquad - (Y_1(2) - 1 - 2Y_2(1) - Y_2(1)^2 + Y_2(2))(Y_1(1)+1+Y_2(1)) + \dfrac{1}{3}(Y_1(1)+1+Y_2(1))^3\\
&\vdots
\end{aligned}\right\}, \qquad (30)$$

$$\left.\begin{aligned}
&F_2(0) = 2,\ F_2(1) = -Y_1(1) - Y_2(1) - 1\\
&F_2(2) = \dfrac{1}{2}Y_1(1)^2 - Y_1(2) - Y_2(2) + \dfrac{1}{2}(1+Y_2(1))^2\\
&F_2(3) = -\dfrac{1}{4}Y_1(1)^3 + Y_1(2)Y_1(1) - Y_1(3) - Y_2(3) + (1+Y_2(1))Y_2(2) - \dfrac{1}{3}(1+Y_2(1))^3\\
&\vdots
\end{aligned}\right\}, \qquad (31)$$

Hence, the series solutions are obtained and given as

$$\tilde{y}_1(t) = 2 - 2t + t^2 - \dfrac{1}{3}t^3 + \dfrac{1}{12}t^4 - \dfrac{1}{60}t^5 + \ldots$$

$$\tilde{y}_2(t) = 1 + \frac{1}{2}t^2 + \frac{1}{6}t^3 + \frac{1}{24}t^4 + \frac{1}{120}t^5 + \ldots$$

For sufficiently large number of terms, the closed forms of the solutions are $\tilde{y}_1(t) = 2e^{-t}$, $\tilde{y}_2(t) = e^t - t$, which are the exact solutions. Table 6 shows the absolute error obtained for three various numbers of terms and at some test points.

**Table 6:** Numerical results for Example 7

| $t_i$ | $|y_1(t_i) - \tilde{y}_1(t_i)|$ | | | $|y_2(t_i) - \tilde{y}_2(t_i)|$ | | |
|---|---|---|---|---|---|---|
| | $N = 5$ | $N = 10$ | $N = 15$ | $N = 5$ | $N = 10$ | $N = 15$ |
| 0.0 | 0.00000000 | 0.00000000 | 0.00000000 | 0.00000000 | 0.00000000 | 0.00000000 |
| 0.2 | 1.7282e-007 | 1.1102e-015 | 1.0521e-016 | 9.1494e-008 | 6.6613e-016 | 2.2204e-016 |
| 0.4 | 1.0759e-005 | 2.0333e-012 | 4.4409e-016 | 6.0310e-006 | 1.0871e-012 | 2.2204e-016 |
| 0.6 | 1.1927e-004 | 1.7309e-010 | 2.2204e-016 | 7.0800e-005 | 9.5652e-011 | ٤,٥١٠١e-016 |
| 0.8 | 6.5259e-004 | 4.0337e-009 | 2.5535e-015 | 4.1026e-004 | 2.3048e-009 | 1.9984e-015 |
| 1.0 | 2.4255e-003 | 4.6229e-008 | 9.0039e-014 | 1.6152e-003 | 2.7313e-008 | 5.0848e-014 |

Tables 2–5 show the absolute errors of the present method for different approximation order, $N$, of the solution and at some test points $t_i$. The results show that the obtained solutions are accurate and converge to the exact ones with increasing the order $N$.

## 5. Conclusions

In this paper, a new general formula and hence new recurrence relations have been derived for computing the differential transform of any analytic nonlinear non-autonomous functions with one or multi-variable. As a special case of the present study, i.e., for autonomous function, the current formulas and recurrence relations reduces to the same mathematical structure as the Adomian polynomials but with constants instead of variable components in the literature. It was found that, the suggested method deals directly with the nonlinear non-autonomous functions, where special transformation or algebraic manipulations were completely avoided. In addition, the method has been successfully applied on different types of differential and intego-differential equations. Moreover, the obtained series solutions demonstrate the validity and applicability of the present approach. Numerically, the results showed that a fast convergence has been achieved for the obtained series solutions. An advantage of the present method is that it can be combined with Padé Approximant (PA), Aftertreatment Techniques (AT), Power Series Extender Method (PSEM) and multi-step technique, among many others. Finally, the authors believe that the present study should be extended to include similar differential and intego-differential equations in the applied sciences, which increases its applicability.

**Conflicts of Interest:** Author has declared that no competing interests exist.